\documentclass[fleqn, 11pt]{frank}

\usepackage{amsmath, amscd, amsthm, amsfonts, amssymb}
\usepackage{mathrsfs, bbm, cancel, stmaryrd, array, comment} 
\usepackage[FIGTOPCAP]{subfigure}
\usepackage[all]{xy}
\usepackage{graphicx,enumitem}
\usepackage{scrtime} 
\usepackage{rotating, lscape, wrapfig, tabularx, caption, color, microtype}
\usepackage[hyperfootnotes=false]{hyperref}
\usepackage{marvosym} 
\usepackage{color}
\usepackage{media9}
\definecolor{darkblue}{rgb}{0,0,1}
\hypersetup{pdfpagemode=UseNone, colorlinks=true, breaklinks=true, linkcolor=blue, menucolor=black, urlcolor=blue, citecolor=red}
\usepackage[bottom]{footmisc}
\usepackage{datetime}

\newcommand{\RR}{\mathbbm{R}}
\newcommand{\CC}{\mathbbm{C}}
\newcommand{\NN}{\mathbbm{N}}

\newtheorem{proposition}{Proposition}
\newtheorem{theorem}[proposition]{Theorem}
\newtheorem{lemma}[proposition]{Lemma}

\newtheorem{corollary}[proposition]{Corollary}

\theoremstyle{remark}

\theoremstyle{definition}

\newtheorem{remark}[proposition]{Remark}
\newtheorem{example}[proposition]{Example}

\allowdisplaybreaks

\begin{document}

\sloppy \raggedbottom

\title{A note on the regularity of matrices with uniform polynomial entries}

\begin{start}

	\author{Frank Klinker{\renewcommand{\thefootnote}{\fnsymbol{footnote}}
				\footnote{Corresponding author}}}{1},\ 
	\author{Christoph Reineke}{2}\\

	\address{Faculty of Mathematics, TU Dortmund University, 44221 Dortmund, Germany\\[0.5ex]
		\href{mailto:frank.klinker@math.tu-dortmund.de}{frank.klinker@math.tu-dortmund.de}\\}{1}
	\address{Im Grubenfeld 12, 44135 Dortmund, Germany\\[0.5ex]
		\href{mailto:christoph_reineke@gmx.de}{christoph\_reineke@gmx.de}\\}{2}

\noindent
\makebox[0.8\textwidth]{%
\begin{minipage}{0.85\textwidth}
\begin{Abstract}
	In this text we study the regularity of  matrices with special polynomial entries. 
Barring some mild conditions we show that these matrices are regular if a natural limit size is not exceeded.
The proof draws connections to generalized Vandermonde matrices and Schur polynomials that are discussed in detail.
	\end{Abstract}
\end{minipage}}

\renewcommand{\dateseparator}{-}
\let\thefootnote\relax



\footnote{\hspace*{-0.5em}Published in {\em S\~ao Paulo Journal of Mathematical Sciences},  \href{https://doi.org/10.1007/s40863-017-0084-6}{doi.org/10.1007/s40863-017-0084-6}}

\end{start}	
\runningheads{%
{\small\sf Preprint} \hspace*{4cm} 
On the regularity of matrices with uniform entries
}{%
{\small\sf Preprint}\hspace*{8.8cm}
F.~Klinker, C.~Reineke
}

\newcommand{\bfx}{\mathbf{x}}
\newcommand{\bfy}{\mathbf{y}}
\newcommand{\bfr}{\mathbf{r}}

\section{Introduction}\label{sec:0}

The content of this text is motivated by an observation on a special type of integer valued matrices: Consider a quadratic matrix that is filled row-wise from the upper left entry to the lower right entry with consecutive integers to a fixed power. Then you see that ''small'' matrices are regular and ''large'' matrices are singular. See Example \ref{ex:1} where the limit value from ''small'' to ''large'' is discussed as an application of Theorem \ref{satz:1}.

As we will see, such matrices appear as a special case of the following general  situation: We consider a non negative integer $\ell$ and three complex valued sequences $\mathbf{x}=(x_1,x_2,\ldots), \mathbf{y}=(y_1,y_2,\ldots)$, and	$\mathbf{r}=(r_1,r_2,\ldots)$ where we assume $\bfr$ to be injektive, i.e.\ $r_i\neq r_j$ for all $i,j\in\NN$.
For each $k\in\NN$ these data define a complex $k\times k$-matrix by restricting the sequences to their first $k$ entries:
\begin{equation}\label{eq:main1}
A=A(k;\bfx,\bfy ,\bfr,\ell)=\begin{pmatrix}
(x_1+r_1y_1)^\ell & (x_1+r_2y_1)^\ell  &\cdots & (x_1+r_ky_1)^\ell  \\
(x_2+r_1y_2)^\ell & (x_2+r_2y_2)^\ell  &\cdots & (x_2+r_ky_2)^\ell  \\
\vdots & \vdots &&\vdots \\
(x_k+r_1y_k)^\ell & (x_k+r_2y_k)^\ell  &\cdots & (x_k+r_ky_k)^\ell  
\end{pmatrix}
\end{equation}
i.e.\
\begin{equation*}\label{eq:main2}
A_{ij}=(x_i+r_jy_i)^\ell\,.
\end{equation*}

The aim of this note is to discuss the regularity of matrices of the form \eqref{eq:main1}. When we do this we will see that there is a strong relation to generalized Vandermonde matrices and we will use a connection to Schur polynomials to get a general result, see Theorem \ref{prop:2}.

We will shortly note two aspects where matrices of the above type or similar occur. 
In general, matrices depending on sequences have applications in combinatorics, graph theory and optimization, see \cite{Brunat,Brunat2} for example.

A more special application of the matrix that we consider is the following: Each row of \eqref{eq:main1}  can be interpreted as a family of ridge functions. This is a family of multivariate functions
\[
\sigma^i(\bfy):=f^i(\mathbf{a}^{i}\cdot\bfy)
\]
depending on fixed sequences $\mathbf{a}^{i}$ and where $\mathbf{a}^{i}\cdot\bfy= \sum_\ell a^i_\ell y_\ell$ denotes the scalar product.\footnote{Usually, if we deal with infinite sequences we have to include some convergence conditions like $\ell^1$, i.e.\ $\sum_k|a_k|<\infty$. For the particular case we are interested in the sums are finite. } In our situation each column of the matrix family of ridge functions of the form
\[
\sigma^{ij}(\bfy)=f^{ij}(\mathbf{a}^{ij}\cdot\bfy)=(x_i+\mathbf{a}^{ij}\cdot\bfy)^\ell
\]
with $\mathbf{a}^{ij}=(\mathbf{e}_j\cdot \bfr)\mathbf{e}_i$. Here we denote  by $\mathbf{e}_j = (\delta_{ij})_{i\in\NN}$ the sequence with all entries vanishing but the $j$-th. Ridge functions have applications in the theory of interpolation and approximation, see \cite{Braess,Pinkus} for more details on this topic.

\section{The general matrix and first remarks}\label{sec:1}

We consider a matrix of the form \eqref{eq:main1} and write $A=A(k;\bfx,\bfy,\bfr,\ell)$ with $A_{ij}=(x_i+r_jy_i)^\ell$. 
In this form the entries of $A$ can be seen as being obtained from polynomials $P_1,\ldots,P_k\in\CC[x]$ of degree at most $\ell$ in the form 
\[
A_{ij}=P_i(r_j)\,.
\]
For such matrices the test on regularity can be very hard, see the nice texts \cite{Kratt,Kratt2} for example. Usually there are connections to representation theory and the theory of planar partitions and we will give further details on this fact in the next section. 
However, in our situation \eqref{eq:main1} one part of the question on regularity is easy to answer:  
The dimension of the space of polynomials of degree at most $\ell$ is $\ell+1$. Therefore, if $k>\ell+1$ there exist complex numbers $\beta_1,\ldots,\beta_k$ with $(\beta_1,\ldots,\beta_k)\neq(0,\ldots,0)$ with $\sum_{i=1}^k\beta_iP_i = 0$. In terms of the rows $R_i(A)$  of $A$ this is $\sum_{i=1}^k\beta_iR_i(A)=0$ such that the rank of $A$ is less than $k$.
The result is summarized in the following Theorem \ref{satz:1}.

\begin{theorem}\label{satz:1}
Independent of the complex valued sequences $\bfx,\bfy,\bfr\in\CC^\NN$ with $\bfr$ injective the matrix  $A(k;\bfx,\bfy,\bfr,\ell)\in M_k(\CC)$ from \eqref{eq:main1} is singular for $k\geq\ell+2$.
\end{theorem}

For the particular case $r_j=j$ we will be more specific here: We will manipulate the columns of $A$ in such a way, that the resulting matrix, $A'$ say, admits at least two equal columns.
 
Therefore, we consider the following operation respecting the determinant: For $j=k,k-1,\ldots,\ell+1$ we replace the $j$-th column $C_j(A)$ by 
\begin{align*}
C_j(A')&
=(-1)^\ell \sum_{\nu=0}^\ell (-1)^{\nu} { \textstyle\binom{\ell}{\nu}} C_{j-\ell+\nu}(A)\,.
\end{align*}
We calculate the $i$-th component of $C_j(A')$ and get 
\begin{align}
 (C_j(A'))_i
=\ & (-1)^\ell\sum_{\nu=0}^\ell (-1)^\nu {\textstyle \binom{\ell}{\nu}} \big(x_i+(j-\ell+\nu)y_i\big)^\ell\nonumber\\
=\ &
\begin{cases}
\ \displaystyle (-1)^\ell y_i^\ell\sum_{\nu=0}^\ell (-1)^\nu {\textstyle \binom{\ell}{\nu}} (a_{ij}+\nu)^\ell & \text{ if }y_i\neq 0\\[1.5ex]
\ \displaystyle (-1)^\ell x_i^\ell\sum_{\nu=0}^\ell (-1)^\nu {\textstyle \binom{\ell}{\nu}}  & \text{ if }y_i= 0
\end{cases}
 \label{eq:2}
\end{align}
with the abbreviation $ a_{ij}:=\frac{x_i}{y_i}+(j-\ell)$ independent of $\nu$. 
For the first case in \eqref{eq:2} the binomial identity yields 
\begin{align*}
(C_j(A'))_i
=\ &(-1)^\ell y_i^\ell\sum_{\mu=0}^\ell{\textstyle\binom{\ell}{\mu}}a_{ij}^{\ell-\mu}\left(\sum_{\nu=0}^\ell(-1)^\nu{\textstyle\binom{\ell}{\nu}}\nu^\mu\right)\,.
\end{align*}
We use \eqref{eq:1} from Lemma \ref{lem:1} and get 
\begin{align}\label{eq:eqcol}
 C_j(A')
= \begin{pmatrix}\ell!\, y_1^\ell\\\vdots\\\ell!\,y_k^\ell\end{pmatrix}\
\end{align}
in both cases of \eqref{eq:2}. 
Therefore, all changed columns $C_j(A')$ coincide for $j=\ell+1,\ldots, k$. 
For at least two columns to coincide we need $k>\ell+1$ which is again the condition from  Theorem \ref{satz:1}.

An even more special situation covers our introductory example:
\begin{example}\label{ex:1}
Let $N$ be a natural number. The sequences $\bfx,\bfy,\bfr$ with  
\[
x_i=N+(i-1)k-1\,,\  y_i=1\,,\ \text{and }\ r_i=i
\] 
yield
\begin{equation}\label{eq:3}
A=
 \begin{pmatrix}
N^\ell & \cdots & (N+k-1)^\ell\\
 (N+k)^\ell &  \cdots & (N+2k-1)^\ell\\
\vdots & &\vdots \\
 (N+(k-1)k)^\ell &  \cdots & (N+k^2-1)^\ell
\end{pmatrix}\,
\end{equation}
which is singular for $k>\ell+1$. 
Its entries are row wise given by the $\ell$-th power of successive natural numbers starting from $N$.
Because $\bfy$ is constant all components of the equal columns $C_j(A')$ in \eqref{eq:eqcol} coincide, too. For $\ell=1$ the matrix  \eqref{eq:3} is of constant-gap type, see \cite{constgap}, and the manipulated matrix $A'$ is
\begin{equation*}\label{eq:3''}
A'= \begin{pmatrix} 
N & 1&\cdots & 1\\
N+k & 1 & \cdots & 1\\
\vdots & \vdots &&\vdots \\
N+(k-1)k & 1 & \cdots & 1
\end{pmatrix}\,.
\end{equation*}
\end{example}

\section{Regularity and relations to representation theory}\label{sec:2}

In this section we will address the following question:
Is the condition on the size of the matrix \eqref{eq:main1} from Theorem \ref{satz:1} sharp? I.e.\ is the matrix $A(k;\bfx,\bfy,\bfr,\ell)$ regular if $k\leq \ell+1$?
For $\ell=0$ this is obviously true. 
Of course, the answer to the question above is ''No'' in general: For $\ell>0$ the trivial examples $\bfx \equiv0$ or $\bfy\equiv 0$ or the example of two sequences  with $x_iy_j=x_jy_i$ for some $1\leq i,j\leq k$ with $i\neq j$ are some counterexamples. 
Nevertheless, the question is reasonable for the matrix \eqref{eq:3} from Example \ref{ex:1}. 

We will follow here a brute-force approach: We will calculate the determinant of $A$. For the moment we will  assume $x_iy_j\neq x_jy_i$ for all $i\neq j$ in \eqref{eq:main1} and the following abbreviations turn out to be useful:
\begin{align*}\begin{aligned}
A_{ij}=\ & 
A_{ij}^{(0)} = x_i^\ell+r_j A_{ij}^{(1)}= b^{(0)}_i+r_j A_{ij}^{(1)}\\
A_{ij}^{(1)}=\ & \binom{\ell}{1}x_i^{\ell-1}y_i +r_jA_{ij}^{(2)} =b^{(1)}_i+r_jA_{ij}^{(2)} \\
 & \vdots\\
A_{ij}^{(\ell-1)}=\ & \binom{\ell}{\ell-1}x_iy_i^{\ell-1} +r_jA_{ij}^{(\ell)}
=b^{(\ell-1)}_i+r_jA_{ij}^{(\ell)}\\
A_{ij}^{(\ell)}=\ & y_i^\ell =b^{(\ell)}_i\,.
\end{aligned}\end{align*}
In the first case $k>l+1$ we will again recover the result from Theorem \ref{satz:1}. 
We use the multilinearity of the determinant and get\footnote{In this compact notation $i$ and $j$ will always number the rows and columns of the matrix, respectively.} 
\begin{align*}
\det(A)=\ &
\begin{vmatrix} 
 \big(x_i^\ell+r_1 A_{i1}^{(1)}\big)_{\substack{i=1,\ldots,k}} &
\big(A_{ij}^{(0)}\big)_{\substack{i=1,\ldots,k\\j=2,\ldots,k}}
\end{vmatrix}
\\
\stackrel{}{=}\ &
\begin{vmatrix} 
 \big(b_i^{(0)}\big)_{\substack{i=1,\ldots,k}} &
\big(A_{ij}^{(0)}\big)_{\substack{i=1,\ldots,k\\j=2,\ldots,k}}
\end{vmatrix}
+
r_1\begin{vmatrix} 
 \big( A_{i1}^{(1)}\big)_{\substack{i=1,\ldots,k}} &
\big(A_{ij}^{(0)}\big)_{\substack{i=1,\ldots,k\\j=2,\ldots,k}}
\end{vmatrix}
\\
\stackrel{}{=}\ &
\begin{vmatrix} 
 \big(b_i^{(0)}\big)_{\substack{i=1,\ldots,k}} &
\big(A_{ij}^{(0)}\big)_{\substack{i=1,\ldots,k\\j=2,\ldots,k}}
\end{vmatrix}
+
r_1 \begin{vmatrix} 
 \big(b_i^{(1)}\big)_{\substack{i=1,\ldots,k}} &
\big(A_{ij}^{(0)}\big)_{\substack{i=1,\ldots,k\\j=2,\ldots,k}}
\end{vmatrix}
\\
&\ +
r_1^{2}\begin{vmatrix} 
 \big( A_{i1}^{(2)}\big)_{\substack{i=1,\ldots,k}} &
\big(A_{ij}^{(0)}\big)_{\substack{i=1,\ldots,k\\j=2,\ldots,k}}
\end{vmatrix}
\\
\vdots\ \ & \\
\stackrel{}{=}\ &
\sum_{\alpha=0}^\ell 
r_1^{\alpha}\begin{vmatrix} 
\big(b_i^{(\alpha)}\big)_{\substack{i=1,\ldots,k}} &
\big(A_{ij}^{(0)}\big)_{\substack{i=1,\ldots,k\\j=2,\ldots,k}}
\end{vmatrix}
\\
\vdots\ \ &  \\
\stackrel{}{=}\ &
\sum_{\alpha_1,\ldots,\alpha_{\ell+1}=0}^\ell 
r_1^{\alpha_1}\cdots r_{\ell+1}^{\alpha_{\ell+1}} \cdot\\ 
&\qquad\qquad\cdot
\begin{vmatrix} 
\big(b_i^{(\alpha_1)}\big)_{\substack{i=1,\ldots,k}} \ldots  \big(b_i^{(\alpha_\ell)} \big)_{\substack{i=1,\ldots,k}}
&\big(A_{ij}^{(0)}\big)_{\substack{i=1,\ldots,k\\j=\ell+2,\ldots,k}}
\end{vmatrix}\,.
\end{align*}
Using the skew symmetry of the determinant this can be written as
\begin{multline*}
\det(A)=
V_{\ell+1}\big(r_1,\ldots,r_{\ell+1}\big) 
\begin{vmatrix} 
\big(b_i^{(0)}\big)_{\substack{i=1,\ldots,k}} \ldots  \big(b_i^{(\ell)} \big)_{\substack{i=1,\ldots,k}}
\big(A_{ij}^{(0)}\big)_{\substack{i=1,\ldots,k\\j=\ell+2,\ldots,k}}
\end{vmatrix}
\end{multline*}
with 
\begin{equation}\label{vand}\begin{aligned}
V_m(u_1,\ldots,u_m)
& =\sum_{\sigma\in {\rm S}_{m}}(-1)^\sigma u_1^{\sigma(1)}\cdots u_m^{\sigma(m)}\\
& = \prod_{1\leq i<j\leq m}(u_j-u_i) 
  =\det\big(u_i^j\big)_{\substack{i=1,\ldots,m\\j=0,\ldots,m-1}}
\end{aligned}
\end{equation}
being the Vandermonde determinant, see \cite{Lorenz} and the discussion below.
Because $k>\ell+1$ we can calculate one more step and get 
\begin{multline*}
\det(A)=V_{\ell+1}\big(r_1,\ldots,r_{\ell+1}\big)
 \sum_{\alpha=0}^\ell r_{\ell+2}^{\alpha}\cdot\\
\cdot\begin{vmatrix} 
\big(b_i^{(0)}\big)_{\substack{i=1,\ldots,k}}& \ldots & \big(b_i^{(\ell)} \big)_{\substack{i=1,\ldots,k}}&\big(b_i^{(\alpha)}\big)_{i=1,\ldots,k}
&\big(A_{ij}^{(0)}\big)_{\substack{i=1,\ldots,k\\j=\ell+3,\ldots,k}}
\end{vmatrix}=0
\end{multline*}
because in each summand at least two columns coincide. 

In the case $k\leq \ell+1$ we can do the same calculations but we will end up with
\begin{align}
\det(A)&= 
\sum_{\alpha_1,\ldots,\alpha_{k}=0}^\ell  \!\!\!
r_1^{\alpha_1}\cdots r_k^{\alpha_{k}}
\begin{vmatrix} 
\big(b_i^{(\alpha_1)}\big)_{\substack{i=1,\ldots,k}} \ldots  \big(b_i^{(\alpha_k)} \big)_{\substack{i=1,\ldots,k}}
\end{vmatrix}\nonumber\\
&=:\sum_{0\leq \alpha_1<\ldots<\alpha_{k}\leq\ell}  \bigg(\sum_{\sigma\in {\rm S}_k} (-1)^\sigma
r_1^{\sigma(\alpha_1)}\cdots r_k^{\sigma(\alpha_{k})}\bigg) \det(B_\alpha)\,.
\label{eq:detA}
\end{align}
Each matrix of which we consider the determinants on the right hand side of \eqref{eq:detA} is a submatrix of the $k\times(\ell+1)$-matrix
\begin{equation}\label{eq:B-basic}
\big(b_i^{(j)}\big)_{\substack{i=1,\ldots,k\\j=0,\ldots,\ell}}
=\begin{pmatrix}
x_1^\ell & \binom{\ell}{1}x_1^{\ell-1}y_1 & \cdots & \binom{\ell}{\ell-1}x_1y_1^{\ell-1} & y_1^\ell\\
\vdots &\vdots &&\vdots &\vdots&\\
x_k^\ell & \binom{\ell}{1}x_k^{\ell-1}y_k & \cdots & \binom{\ell}{\ell-1}x_ky_k^{\ell-1} & y_k^\ell
\end{pmatrix}\,.
\end{equation}
It is defined by a sequence $\alpha=(\alpha_1,\ldots,\alpha_k)$ of strictly increasing natural numbers $0\leq \alpha_1< \ldots <\alpha_k\leq \ell$ and given by
\begin{align*}
B_\alpha=\big(b_i^{(\alpha_j)}\big)_{\substack{i,j=1,\ldots,k}}
=\begin{pmatrix}
\binom{\ell}{\alpha_1}x_1^{\ell-\alpha_1} y_1^{\alpha_1} &\ldots & \binom{\ell}{\alpha_k}x_1^{\ell-\alpha_k} y_1^{\alpha_k} \\
\vdots & &\vdots&\\
\binom{\ell}{\alpha_1}x_k^{\ell-\alpha_1} y_k^{\alpha_1} &\ldots & \binom{\ell}{\alpha_k}x_k^{\ell-\alpha_k} y_k^{\alpha_k} 
\end{pmatrix}\,.
\end{align*}
Its determinant is given by
\begin{align}
\det(B_\alpha)&=
\prod_{i=1}^k x_i^\ell \prod_{j=1}^k\binom{\ell}{\alpha_j}
\begin{vmatrix}
\rho_1^{\alpha_1} & \ldots &\rho_1^{\alpha_k}\\
\vdots&&\vdots\\
\rho_k^{\alpha_1} & \ldots &\rho_k^{\alpha_k}
\end{vmatrix}\nonumber\\
&=:\prod_{i=1}^k x_i^\ell \prod_{j=1}^k\binom{\ell}{\alpha_j}
V_{k,\alpha}(\rho_1,\ldots,\rho_k)
\label{eq:detBa}
\end{align}
with $\rho_i:=\frac{y_i}{x_i}$ for all $1\leq i\leq k$. Interchanging the role of $\bfx$ and $\bfy$ yields the following formula for $\det (B_\alpha)$ where we consider the strictly increasing sequence $\alpha^\complement=(\ell-\alpha_k,\ldots,\ell-\alpha_1)$ associated to $\alpha$:
\begin{align}
\det(B_\alpha)
&=\prod_{i=1}^k y_i^{-\ell} \prod_{j=1}^k\binom{\ell}{\alpha_j}
\begin{vmatrix}
(\frac{1}{\rho_1})^{\ell-\alpha_1} & \ldots &(\frac{1}{\rho_1})^{\ell-\alpha_k}\\
\vdots&&\vdots\\
(\frac{1}{\rho_k})^{\ell-\alpha_1} & \ldots &(\frac{1}{\rho_k})^{\ell-\alpha_k}
\end{vmatrix}\nonumber\\
&=(-1)^{\frac{k(k-1)}{2}}\prod_{i=1}^k y_i^{-\ell} \prod_{j=1}^k\binom{\ell}{\alpha_j}
V_{k,\alpha^\complement}\big(\tfrac{1}{\rho_1},\ldots,\tfrac{1}{\rho_k}\big)\,.
\tag{\ref{eq:detBa}'}\label{eq:detBa'}
\end{align}
Moreover, the factor in \eqref{eq:detA} we put in round brackets is also of this special form, and we may write by using  \eqref{eq:detBa} or \eqref{eq:detBa'}
\begin{equation}
\det(A)  = \prod_{i=1}^k x_i^{\ell}  \sum_{0\leq \alpha_1<\ldots<\alpha_{k}\leq\ell}\, \prod_{j=1}^k{\binom{\ell}{\alpha_j}} 
	 V_{k,\alpha}(r_1,\ldots,r_k)\, V_{k,\alpha}\big(\rho_1,\ldots,\rho_k\big)\,\label{detAneu}
\end{equation}
or
\begin{multline}
\det(A)  = (-1)^{\frac{k(k-1)}{2}}\prod_{i=1}^k {y_i^{\ell}}\ \cdot \\ \cdot\sum_{0\leq \alpha_1<\ldots<\alpha_{k}\leq\ell}\, \prod_{j=1}^k{\binom{\ell}{\alpha_j}} 	V_{k,\alpha}(r_1,\ldots,r_k)\, V_{k,\alpha^\complement}\big(\tfrac{1}{\rho_1},\ldots,\tfrac{1}{\rho_k}\big)\,.\tag{\ref{detAneu}'}\label{detAneu2}
\end{multline}
The remaining determinant $V_{k,\alpha}(u_1,\ldots,u_k)$ is a generalized Vandermonde determinant, see \cite{Heinemann}. 
In particular, $V_{k,(0,1,\ldots,k-1)}=V_k$. Taking into account that $V_{k,\alpha}$ and $V_k$ have common zeros, their quotient is a polynomial, too. All these polynomials are symmetric and given by the Schur polynomials, see e.g.~\cite{DeMarchi,King,Littlewood}. 
If we associate to the strictly increasing sequence $\alpha=(\alpha_1,\ldots,\alpha_k)$ the non increasing sequence $\lambda=(\lambda_1,\ldots,\lambda_k)$ with $\lambda_i:=\alpha_{k-i+1}-k+i$ the mentioned Schur polynomial is given by the polynomial $s_{\lambda}(x_1,\ldots, x_k)$ with 
\begin{equation*}
s_{\lambda}(u_1,\ldots,u_k) V_k(u_1,\ldots,u_k)  = V_{k,\alpha}(u_1,\ldots,u_k)=\det \big( (u_i^{\alpha_j})_{i,j=1,\ldots,k}\big)\,.
\end{equation*}
In particular, they are used to evaluate the character and -- more applied -- the dimension of representations of the classical Groups $GL(n)$. In this notation the latter are described by Young tableaux of shape $\lambda$, see \cite{FultonHarris,ChipMoh,KoikeTerada,Koike} for details on this topic or \cite{Klinker} for an application. There is one important fact that we would like to recall here -- for a discussion of this fact we like to refer to the nice pair of papers \cite{Stanley}:
\begin{theorem}
Let $s_\lambda$ be the Schur  polynomial associated to the non increasing sequence $\ell\geq \lambda_1\geq\lambda_2\geq\ldots\geq\lambda_k\geq 0$. Then the expansion of $s_\lambda$ with respect to monomials is of the form 
\begin{equation*}
s_\lambda(u_1,\ldots,u_k)=\sum_\mu \Gamma_\lambda^\mu \ u_1^{\mu_1}\cdots u_k^{\mu_k}
\end{equation*}
where the sum is taken over all $k$-tuples $\mu=(\mu_1,\ldots,\mu_k)$. The coefficient $\Gamma_\lambda^\mu$ is obtained as follows: 
Identify $\lambda$ with the Young tableau with $\lambda_i$ boxes in the $i$-th row. Now fill the boxes with $\mu_1$ times '1', $\mu_2$ times '2',$\ldots$, $\mu_k$ times '$k$' such that the rows are non decreasing and the columns are strictly increasing. The coefficient $\Gamma_\lambda^\mu$ is now given by the number of ways this can be done. 
\end{theorem}
\begin{corollary}\label{cor:1}
All coefficients in the expansion of Schur polynomials with respect to monomials are non negative. Therefore, the Schur polynomials are positive if we restrict  to positive real values.
\end{corollary}
\begin{corollary}\label{cor:2}
Suppose $\rho_i=\frac{y_i}{x_i}\in\RR^+$. Then the determinants of the matrices $B_\alpha$ from \eqref{eq:detBa}  do not vanish and so does not $\det(A)$ from \eqref{detAneu} if we assume $r_i>0$ in addition.
\end{corollary}

We will shortly discuss the cases $x_i=0$ or $y_i=0$ for one $1\leq i\leq k$. First we consider the case $x_1=0$.
Then \eqref{eq:B-basic} yields  
$\det(B_\alpha)=0$ whenever $\alpha_k\neq \ell$. Therefore, if we write $\alpha=(\alpha',\ell)$ with $\alpha'=(\alpha_1,\ldots,\alpha_{k-1})$ we get
\begin{equation}
\det (B_{(\alpha',\ell)})=(-1)^k y_1^\ell \prod_{i=2}^{k} x_i^\ell \prod_{j=1}^{k-1}{\binom{\ell}{\alpha_j}} V_{k-1,\alpha'}(\rho_1,\ldots,\rho_{k-1})\,.\label{16}
\end{equation}
In the same way for $y_i=0$ we get $\det(B_\alpha)=0$ whenever $\alpha_1\neq 0$. Therefore, if we write $\alpha=(0,\alpha')$ with $\alpha'=(\alpha_2,\ldots,\alpha_k)$ we get 
\begin{equation}\label{16'}
\det (B_{(0,\alpha')})= x_1^\ell \prod_{i=2}^k x_i^\ell \prod_{j=2}^k{\binom{\ell}{\alpha_j}}
V_{k-1,\alpha'}(\rho_2,\ldots,\rho_{k})\,. \tag{\ref{16}'}
\end{equation}
In both cases we end up with the same type of determinant but of size one less. Therefore, Corollary \ref{cor:2} holds in these cases, too.

\begin{example}[Example \ref{ex:1} continued]\label{ex:cont}
We recall that in this example $x_i=N+(i-1)k-1$, $y_i=1$ and $r_i=i$. 
For the regularity question applied to matrix \eqref{eq:3} we have to distinguish the two cases $N>1$ and $N=1$  corresponding to $x_i\neq 0$ for all $1\leq i\leq k$ and $x_1=0$, respectively. Corollary \ref{cor:2} and the remark thereafter show that $A$ from \eqref{eq:3} is regular for $k\leq\ell+1$. The determinant is a polynomial of degree ${\ell k}$ in $N$, a priori. However, in the limiting case $k=\ell+1$ we obtain degree 0 such that the determinant is independent of $N$. 
This is due to the fact that in this case only one summand in \eqref{detAneu2} is left, namely the one with $\alpha=\alpha^\complement=(0,1,\ldots\ell)$. If we expand the Vandermonde determinants as in \eqref{vand} the latter is given by
\begin{align*}
\det(A) 
&=(-1)^{\frac{\ell(\ell+1)}{2}} (\ell+1)^{\frac{\ell(\ell+1)}{2}}
	\prod_{j=0}^\ell (j!)^2{\binom{\ell}{j}}\,.
\end{align*}
\end{example}

We summarize the discussion in the following Theorem \ref{prop:2} which, in particular, covers Examples \ref{ex:1} and \ref{ex:cont}.
\begin{theorem}\label{prop:2}
Let $\bfx,\bfy$ be complex valued sequences with $\frac{x_i}{y_i}\in\RR^+$ and $x_iy_j-y_ix_j\neq 0$ for all $1\leq i,j\leq k$, and $\bfr$ be an positive injective  sequence.
Then the matrix $A=A(k;\bfx,\bfy,\bfr,\ell)\in M_k\CC$ 
with
\begin{equation*}
A_{ij}=\big(x_i+r_jy_i\big)^\ell 
\end{equation*}
is regular if and only if $k\leq \ell+1$.
\end{theorem}
We end this note by some remarks 
\begin{remark}
\begin{enumerate}
\item 
The condition on $\bfx$ and $\bfy$ from Theorem \ref{prop:2} may be relaxed by letting $x_i$ or $y_i$ vanish for one $i$, see \eqref{16} and \eqref{16'}.
\item 
In fact, Theorem \ref{prop:2} is true for a generic choice of sequences $\bfx, \bfy$. 
The values for which the matrix fails to be regular are the solutions of a polynomial equation and are connected to Vandermonde varieties. More or less, the latter are given by the set of zeros of $V_{k,\alpha}(q_1,\ldots,q_k)$ which are of codimension one. See \cite{3} for some calculations regarding such varieties, also in the case of non quadratic Vandermonde-type matrices.
\item 
The calculations in this note that make use of the relation between the regularity of matrices and determinants are valid in any algebraically closed field and, therefore, we may replace $\CC$ by such field. For an application of Vandermonde determinants over such fields and their connection to linear recurrence sequences see \cite{4}.
\end{enumerate}
\end{remark}

\section{Appendix: Some useful calculations} \label{sec:app}

When we simplified \eqref{eq:2} we used the following observation.
\begin{lemma}\label{lem:1}
For any $\ell\in\NN$ and polynomial $q(x)$ of degree at most $\ell$ we have
\begin{align}\label{eq:1}
\sum_{\nu=0}^\ell(-1)^\nu{\binom{\ell}{\nu}}q(\nu) = (-1)^\ell q^{(\ell)}(0) \,.
\end{align}
\end{lemma}
The result from Lemma \ref{lem:1} is obtained in a way analog to the well known special case $q(x)=1$. We prove the formula by showing that it is true for any monomial up to degree $\ell$. We consider the binomial identity $(1-x)^\ell = \sum_{\nu=0}^\ell(-1)^\nu{\binom{\ell}{\nu}} x^\nu$ and its $\ell$ derivatives 
\begin{align*}
(-1)^{s}\ell\cdot\ldots\cdot(\ell-s+1)\, (1-x)^{\ell-s}
& = \sum_{\nu=0}^\ell(-1)^\nu{\binom{\ell}{\nu}}\nu\cdot\ldots\cdot(\nu-s+1) x^{\nu-s}
\end{align*}
for $s=1,\ldots,\ell$. 
We write  
$
\pi_0(\nu):=1$ and $\pi_s(\nu)=\nu(\nu-1)\cdot\ldots\cdot(\nu-s+1)$ for $s>1$ such that
 $\pi_s$ is a polynomial of degree $s$ with leading coefficient $1$, i.e.\ 
$
\pi_s(\nu)=\sum\limits_{\kappa=0}^{s}b^{(s)}_\kappa\nu^\kappa\,
$ 
for some real coefficients $b^{(s)}_\kappa$ and $b_s^{(s)}=1$. Using this the $\ell+1$ equations above are given by 
\begin{align*}
(-1)^ss!{\binom{\ell}{s}} (1-x)^{\ell-s} 
 = \sum_{\kappa=0}^{s} b_{s-\kappa}^{(s)}\left(\sum_{\nu=0}^\ell(-1)^\nu {\binom{\ell}{\nu}}\nu^{s-\kappa} x^{\nu-s}\right) 
\end{align*}
for $s=0,\ldots,\ell$. 
We insert the special value $x=1$ into all $\ell+1$ equations and get 
\begin{align*}
0 & = \sum_{\nu=0}^\ell(-1)^\nu{\binom{\ell}{\nu}}\,, \quad
0  = \sum_{\nu=0}^\ell(-1)^\nu{\binom{\ell}{\nu}}\nu 
	  	+b^{(1)}_0 \sum_{\nu=0}^\ell(-1)^\nu{\binom{\ell}{\nu}}\,,\quad \ldots\,,\\
(-1)^\ell \ell! &=\sum_{\nu=0}^\ell(-1)^\nu{\binom{\ell}{\nu}}\nu^\ell
		+b^{(\ell)}_{\ell-1}\sum_{\nu=0}^\ell(-1)^\nu{\binom{\ell}{\nu}}\nu^{\ell-1}
+\ldots	+b^{(\ell)}_0\sum_{\nu=0}^\ell(-1)^\nu{\binom{\ell}{\nu}}\,.
\end{align*}
From top to bottom this yields \eqref{eq:1} for the monomials $1,x,x^2,\ldots, x^\ell$ which proves the Lemma.


\begin{thebibliography}{A}


\bibitem{constgap}
Boris Aronov, Tetsuo Asano, Yosuke Kikuchi, Subhas C.\ Nandy, Shinji Sasahara, and Takeaki Uno:
\newblock{A Generalization of Magic Squares with Applications to Digital Halftoning}
\newblock{\em Theory Comput.\ Syst.} {\bf 42} (2008) 143-156

\bibitem{Braess}
Dietrich Braess:
\newblock {\em Nonlinear Approximation Theory} (Springer Series in Computational Mathematics 7).
\newblock Springer-Verlag, 1986

\bibitem{Brunat}
Joseph M.\ Brunat and Antonio Montes:
\newblock{The Power-Composition Determinant and Its Application to Global Optimization.}
\newblock{\em SIAM. J. Matrix Anal. Appl.} {\bf 23} no.~2 (2001) 459–471

\bibitem{Brunat2}
J.\ M.\ Brunat, C.\ Krattenthaler, A.\ Lascoux, and A.\ Montes:
\newblock{Some composition determinants.}
\newblock{\em Linear Algebra Appl.} {\bf 416} (2006) 355-364

\bibitem{ChipMoh}
Jaydeep Chipalkatti and Tagreed Mohammed:
\newblock{Standard Tableaux and Kronecker projections of Specht modules.}
\newblock{\em  Int.\ Electron.\ J.\ Algebra} {\bf 10} (2011) 123-150

\bibitem{DeMarchi}
Stefano De Marchi:
\newblock{Polynomials arising in factoring generalized Vandermonde determinants: an algorithm for computing their coefficients.}
\newblock {\em Math.\ Comput.\ Modelling} {\bf 34} (2001) 271-281

\bibitem{3}
Ralf Fr\"oberg and Boris Shapiro:
\newblock{On Vandermonde varieties.} 
\newblock{\em Math.\ Scand.} {\bf 119} (2016) no.~1 73–91

\bibitem{FultonHarris}
William Fulton and Joe Harris:
\newblock {\em Representation theory} (Graduate Texts in Mathematics, 129).
\newblock Springer-Verlag, New York, 1991

\bibitem{Heinemann}
E.~R.~Heineman: 
\newblock{Generalized Vandermonde determinants.}
\newblock{\em Trans.\ Amer.\ Math.\ Soc.} {\bf 31} no.~3 (1929) 464-476

\bibitem{King}
R.~C.~King:
\newblock{Generalised Vandermonde determinants and Schur functions.}
\newblock{\em Proc.\ Amer.\ Math.\ Soc.} {\bf 48} no.~1 (1975) 53-56

\bibitem{Klinker}
Frank Klinker:
\newblock{The decomposition of the spinor bundle of Grassmann manifolds.}
\newblock{\em J.\ Math.\ Phys.} {\bf 48} (2007) no.~11, 113511, 26 pp

\bibitem{Koike}
Kazuhiko Koike and Itaru Terada: 
\newblock On the decomposition of tensor products of  representations of the classical groups by means of the universal character.
\newblock {\em Adv.\ Math.}, {\bf 74} (1989) 57-86

\bibitem{KoikeTerada}
Kazuhiko Koike and Itaru Terada:
\newblock Young-diagrammatic methods for the representation theory of the
  classical groups of type {$B\sb n,\;C\sb n,\;D\sb n$}.
\newblock {\em J.\ Algebra}, {\bf 107} no.~2 (1987) 466-511

\bibitem{Kratt}
Christian Krattenthaler:
\newblock Advanced determinant calculus.
\newblock{\em S\'emin.\ Lothar.\ Comb.} {\bf 42} (1999) Art.\ B42q, 67 pp (electronic)

\bibitem{Kratt2}
Christian Krattenthaler:
\newblock Advanced determinant calculus: A complement.
\newblock{\em Linear Algebra Appl.} {\bf 411} (2005) 68-166

\bibitem{Littlewood}
D.~E.~Littlewood: 
\newblock{\em The Theory of Group Characters and Matrix Representations of Groups.}
\newblock Oxford University Press, New York, 1940

\bibitem{Lorenz}
Falko Lorenz:
\newblock {\em Lineare Algebra I}. 
\newblock Spektrum Akademischer Verlag, 4.~Aufl.~2003, 2.~Nachdruck 2008

\bibitem{Pinkus}
Allan Pinkus:
\newblock{\em Ridge Functions} (Cambridge Tracts in MAthematics 205).
\newblock Cambridge University Press, 2015

\bibitem{4}
Hans Peter Schlickewei and Carlo Viola:
\newblock{Generalized Vandermonde determinants.}
\newblock{\em Acta Arith.} {\bf 95} (2000) no.~2 123–137

\bibitem{Stanley}
Richard P.~Stanley:
\newblock {Theory and application of plane partitions. Part I} and {Part II}. 
\newblock {\em Studies in Appl.~Math.} {\bf 50} (1971) 167-188 and 259-279

\end{thebibliography}
\end{document}